\documentclass[a4paper,10pt]{article}
\usepackage[utf8]{inputenc}
\usepackage{amssymb}
\newtheorem{theorem}{Theorem}

\newtheorem{lemma}{Lemma}
\newtheorem{example}{Example}
\newtheorem{problem}{Problem}
\title{Small palindromic lengths in free groups and word equations with antimorphisms}
\author{Anna E. Frid\footnote{Aix Marseille Univ, CNRS, I2M, Marseille, France; anna.e.frid@gmail.com}}

\begin{document}

\maketitle

\begin{abstract}
The palindromic length of a finite word $w$ is defined as the minimal number of palindromes such that their product is $w$. Clearly, this function may take different values depending on if we consider $w$ as an element a free semigroup or of a free group: for example, in the free semigroup, the palindromic length of $abca$ is 4 (here every letter is a palindrome), and in the free group, it is 3 since $abca=(aba)(a^{-1}a^{-1})(aca)$. 

In free semigroups, the palindromic length can clearly be computed, and there are fast algorithms for that. In free groups, the question is trickier. In this paper, we characterize words in the free group whose palindromic length is 2 and 3.
\end{abstract}

\section{Introduction}
The mirror image $\overline{X}$ of a word $X$ over an alphabet $\Sigma$ is the word written in the opposite direction: if $X=x_1x_2\cdots x_n$, $x_i \in \Sigma$, then $\overline{X}=x_n\cdots x_2x_1$. A {\it palindrome} is a word equal to its mirror image: $P=\overline{P}$, like $rotator$ or $abba$.

The {\it palindromic length} of a word $X$ is the minimal number of palindromes such that their product is $X$. As it has been discussed in the abstract, this notion differs if we consider $X$ as an element of a free semigroup $\Sigma^*$ or of a free group $\langle\Sigma\rangle$. Clearly, the second case give more ways to decompose a word to palindromes, so, the group palindromic length is not greater than the semigroup one.

The semigroup palindromic length is well-studied since the 2013 paper by Puzynina, Zamboni and the author \cite{fpz}. In particular, it was conjectured in that paper that the palindromic length of factors is unbounded for every aperiodic infinite word; after several partial results, the complete proof of the conjecture has been published by J. Rukavicka \cite{rukavicka}. The semigroup length of a finite word can obviously be computed by listing all possible factorisations; moreover, very fast algorithms have been found for that \cite{rs,bkrs}.

As for the group palindromic length, it was first considered by Bardakov, Shpilrain, Tolstykh in 2005 under the name of palindromic width \cite{bst}. The main result of that paper is that this function is unbounded: the free group palindromic length of words of the form $aba^2b^2\cdots a^nb^n$ grows to infinity with $n$.

At the same time, the difference between the two versions of the palindromic length can be arbitrarily large, as it is shown by Saarela \cite{saarela}: 

\begin{example}[Saarela \cite{saarela}]
Consider an alphabet $\{a_1,\ldots,a_n,b,c\}$ and denote $A=a_1a_2\cdots a_n$. Consider the word
\[w=Abc\overline{A}=a_1a_2\cdots a_n bc a_n \cdots a_2a_1.\]
In the free semigroup, its palindromic length is equal to its length $2n+2$, since there are no palindromes of length more than 1 in it. On the other hand, in the free group, the palindromic length of $w$ is $3$ since
\[w=(Ab\overline{A})(A\overline{A})^{-1}(Ac\overline{A}) .\]
\end{example}

\medskip
Several years ago, V. Bardakov asked in a private communication if the problem of finding the palindromic length is decidable. Since the palindromes can appear ``from nowhere'' as in the example above, this question is not trivial.

The main idea of this paper is to reduce the question if a word is decomposable to $k$  palindromes in a free group to solving word equations in a free semigroup. Unfortunately, since these equations contain antimorphisms and in particular the condition for a word to be a palindrome, we cannot immediately report that they all can be solved with the Makanin algorithm \cite{diekert} or another known method. However, we can solve several relatively simple equations to establish the the main  results of the paper. They are the following characterisations of reduced forms of words with the palindromic complexity in the free group equal to 2 and 3:

\begin{theorem}\label{t:2}
 A word in a free group $\langle\Sigma\rangle$ is a product of at most two palindromes if and only if its reduced form is a concatenation of the form
 \[APQA^{-1}\]
 for some (possibly empty) palindromes $P,Q$ and a word $A$.
\end{theorem}
\begin{theorem}\label{t:3}
 A word in a free group $\langle\Sigma\rangle$ is a product of at most three palindromes if and only if its reduced form is a concatenation of the form
 \[ABPQB^{-1}R \overline{A} \mbox{~or~} APBQRB^{-1}\overline{A},\]
 where $P,Q,R$ are (possibly empty) palindromes and $A,B$ are arbitrary words.
\end{theorem}
Since in both theorems, we characterize the reduced form, not just a group product, as a concatenation of a limited number of words, every time we can check in a polynomial time if a word is of a given form or not. However, even though Theorem \ref{t:2} is relatively simple, the only existing proof of Theorem \ref{t:3} is long, contains dozens of cases and requires solutions of several specific word equations with palindromes and antimorphisms. Analogous results for words of palindromic lengths $k\geq 4$ would have been much too technical, since the number of cases grows at least exponentially.

 So, the following question stays open.
\begin{problem}
 Can we find analogous characterisations of words of any greater palindromic lengths in the free group?
\end{problem}
At the end of the paper, we will discuss possible formulations of such potential results and difficulties with them.

Theorem \ref{t:3} was initially conjectured in a 2019 Master thesis of Moujahid Bou-Laouz; unfortunately, his work was not continued and did not result to a correct proof.

The next long section of the paper is devoted to solving a series of semigroup word equations which will be useful for proving the teorems in Sections \ref{s:t2} and \ref{s:t3}. After that, in a conclusion we make a conjecture on how analogous results for greater $k$ may look like and discuss why the next theorems cannot be similarly obtained by hand.

\section{Word equations with palindromes}
All the results of this section are formulated for elements of a free semigroup, even though they will later be used for reduced forms of words in a free group. To avoid this confusion, we can assume that the words are considered over the alphabet
\[\Sigma_g=\Sigma\cup \{a^{-1}|a\in \Sigma\},\]
and that letters $a$ and $a^{-1}$ are by chance never consecutive. So, the results for the semigroup $\Sigma_g^*$ can directly be applied to reduced forms of words from the free group $\langle \Sigma \rangle$.

An {\it antimorphism} is a map $\varphi:\Sigma^*_g\to\Sigma^*_g$ such that $\varphi(XY)=\varphi(Y)\varphi(X)$ for all words $X$, $Y$. Clearly, an antimorphism is completely defined by images of letters. For the main results of this paper, we will only need two automorphisms in $\Sigma_g^*$: the mirror image $\overline{\;\cdot\;}$ and the inverse defined by $a \mapsto a^{-1}$ and $a^{-1} \mapsto a$ for every symbol $a$ from the initial alphabet $\Sigma$. 

Note also that the condition for a word $P$ to be a palindrome can immediately be formulated in terms of the antimorphism 
$\overline{\;\cdot\;}$: $P=\overline{P}$. Note also that the two bijective antimorphisms considered here commute, and their composition is a bijective morphism: If $X=x_1\cdots x_n$, $x_i \in \Sigma_g$, then $\overline{X}^{-1}=\overline{X^{-1}}=x_1^{-1}\cdots x_n^{-1}$. The same commutation would hold if we consider instead of the inverse any other antimorphism such that $\varphi(\varphi(a))=a$ for all $a$. 

The word equations considered here contain palindromes and the inverse antimorphism, so, they cannot be immediately solved by the classical Makanin algorithm \cite{diekert}, even though we use the same ideas as in the simplest version of it. As for word equations with antimorphisms, the only paper we aware of for now \cite{bmghl} considers equations with prescribed length of words, which is quite far from the general case.


Let us denote the semigroup palindromic length of a word $X$ by $|X|_p$. The following theorem is extremely useful for the proofs below.

\begin{theorem}[Saarela \cite{saarela}]\label{saar}
For all words $X$ and $Y$ over the free semigroup, 
\[|X|_p\leq |Y|_p+|XY|_p \mbox{~and~} |Y|_p\leq |X|_p+|XY|_p.\]
\end{theorem}

In most cases, the Saarela theorem above is used in its simplest form:
\begin{lemma}\label{l:P=QX}
Suppose that
\[P=QX \mbox{~or~} P=XQ,\]
where $P$ and $Q$ are palindromes. Then  $X$ is a concatenation of at most two palindromes.
\end{lemma}

\medskip
The next lemma is a direct corollary of the previous one.
\begin{lemma}\label{l:px=qy}
 Suppose that $PX=QY$, where $P$ and $Q$ ar palindromes. Then either $X=STY$ or $Y=STX$ for some palindromes $S$ and $T$, depending on which of the words $X$ or $Y$ is longer.
\end{lemma}

\medskip
The next lemma is just a simple case study.
\begin{lemma}\label{l:xy=zp}
 Suppose that $XY=ZP$ for some palindrome $P$. Then one of the following holds:
 \begin{enumerate}
  \item either $Z=XA$ and $Y=AP$ for a word $A$;
  \item or $X=ZA$ and $Y=P'\overline{A}$ for a word $A$ and a palindrome $P'$;
  \item or $X=Z\overline{Y}P'$ for a palindrome $P'$.
 \end{enumerate}
\end{lemma}
{\noindent \sc Proof.} This is just a formalization of three possible cases: either $X$ is shorter than $Z$; or longer than $Z$ but does not continue to the middle of $P=AP'\overline{A}$; or goes until the middle of $P=\overline{Y}P'Y$. \hfill $\Box$

The next lemma is just a continuation of Lemma \ref{l:px=qy} involving one factor more.

\begin{lemma}\label{PXY=QZ}
Suppose that $PXY=QZ$, where $P$ and $Q$ are palindromes and $X$, $Y$, $Z$ are arbitrary words. Then one of the following six cases holds: there exist two palindromes $S$ and $T$ and eventually an arbitrary word $X'$ such that
\begin{enumerate}
 \item either $Z=STXY$;
 \item or $X=STX'$ and $Z=X'Y$;
 \item or $Y=S\overline{X}TZ$;
 \item or $X=X'S$ and $Y=\overline{X'}YZ$;
 \item or $X=SX'$ and $Y=T\overline{X'}Z$;
 \item or $X=SX'T$ and $Y=\overline{X'}Z$.
\end{enumerate}\end{lemma}
\noindent {\sc Proof.} From Lemma \ref{l:px=qy} we know that either $Z\in \Sigma_P^2 XY$ (which is exactly the first of six cases), or $XY \in \Sigma_P^2Z$. The last case in its turn is divided into two subcases: either $Y$ is a suffix of $Z$, and then $X\in \Sigma_P^2X'$ and $Z=X'Y$ (it is exactly the second of the six cases), or $Z$ is a suffix of $Y$. If so, we have to consider the four possible positions of the boundary between $X$ and $Y$. 

If $X$ is a prefix of the first of the two palindromes and does not cover its center, we denote the first palindrome by $XS\overline{X}$, the second by $T$ and observe the third case of the lemma.

If $X$ is a prefix of the first of the two palindromes and does cover its center, we denote this first palindrome by $X'S\overline{X'}$ with $X=X'S$ and observe the fourth case of the lemma.

If $X$ is the prefix of the two palindromes containing the first of them but not containing the center of the second one, we denote the first palindrome by $S$, the second by $X'T\overline{X'}$ with $X=SX'$, and get the fifth case of the lemma.

At last, if $X$ is the prefix of the two palindromes containing the first of them and the center of the second one, we again 
denote the first palindrome by $S$, the second by $X'T\overline{X'}$, but this time $X=SX'T$, and so we get the last case. \hfill $\Box$

\bigskip
The next lemma is just an obvious list of cases.
\begin{lemma}\label{l:XY=PQ}
 Suppose that $XY=PQ$ (in reduced forms), where $P$ and $Q$ are palindromes and $X$ and $Y$ are arbitrary words. Then one of the following holds:
 \begin{enumerate}
  \item either $Y=P'\overline{X}Q$, where $P'$ is a center palindrome of $P=XP'\overline{X}$;
  \item or $X=ZP'$ and $Y=\overline{Z}Q$, where $P=ZP'\overline{Z}$;
  \item or $X=PZ$ and $Y=Q'\overline{Z}$, where $Q=ZQ'\overline{Z}$;
  \item or $X=P\overline{Y}Q'$, where $Q=\overline{Y}Q'Y$.
 \end{enumerate}

\end{lemma}

This lemma just uses the previous one.
\begin{lemma}\label{l:PX=YQ}
 Suppose that $PX=YQ$ (in reduced forms), where $P$ and $Q$ are palindromes and $X$ and $Y$ are arbitrary words. Then there exists palindromes $S$, $T$, and eventually and arbitrary word $Z$ such that one of the following holds:
 \begin{enumerate}
  \item either $X=SYT$;
  \item or $X=ZS$ and $Y=TZ$;
  \item or $X=SZ$ and $Y=ZT$;
  \item or $Y=SXT$.
 \end{enumerate}
 \end{lemma}
\noindent {\sc Proof.} Once again, it is just a case study.

Suppose first that $Y$ is longer than $P$, that is, $Y=PZ$. Then $PX=PZQ$ and so $X=ZQ$, falling into the 2nd case.

Now suppose that $Y$ is not longer than $P$ but longer than its half: $Y=ZP'$, where $P=ZP'\overline{Z}$. Then $PX=YQ$ means $\overline{Z}X=Q$, and so there are two subcases: either $Q=\overline{Z}Q'Z$, and then $X=Q'Z$, falling into case 3; or $Q=\overline{X}Q'X$ with $\overline{Z}=\overline{X}Q'$ and thus $Z=Q'X$. So, we have $Y=Q'XP$, falling into case 4.

At last, suppose that $Y$ is shorter than a half of $P$: P=$YP'\overline{Y}$. Then $PX=YQ$ means $P'\overline{Y}X=Q$, and from Lemma \ref{l:P=QX}, $\overline{Y}X$ is a product of two palindromes. Now it remains to apply Lemma \ref{l:XY=PQ} and to return from the mirror image of $Y$ to get exactly the four cases of the statement. \hfill $\Box$

\begin{lemma}\label{l:XPY=Q}
 Suppose that $XPY=Q$ for some palindromes $P$ and $Q$. Then either $X=\overline{Y}ST$ or $Y=ST\overline{X}$ for some palindromes $S$ and $T$.
\end{lemma} \noindent {\sc Proof.} If $X$ is longer than $Y$, then it starts with $\overline{Y}$ since $Q$ is a palindrome: $X=\overline{Y}Z$. Cutting $Y$ and $\overline{Y}$ from the equation, we get $ZP=Q'$, where $Q'$ is a palindrome. By the statement symmetric to Lemma \ref{l:P=QX}, we get $Z=ST$, where $S$ and $T$ are palindromes.

The case of $Y$ longer than $X$ is similar. \hfill $\Box$

\begin{lemma}\label{l:XPY=QZ}
 Suppose that $XPY=QZ$ for some palindromes $P$ and $Q$. Then there exist palindromes $S$ and $T$ and eventually an arbitrary word $F$ such that one of the following holds:
 \begin{enumerate}
 \item either $X=\overline{F}ST$ and $Y=FZ$;
  \item or $Y=ST\overline{X}Z$;
  \item or $Z=SXTY$;
  \item or $X=SF$ and $Z=FTY$;
  \item or $X=FS$ and $Z=TFY$;
  \item or $X=SFT$ and $Z=FY$.
 \end{enumerate}\end{lemma}
\noindent {\sc Proof.} If $Z$ is shorter than $Y$, $Y=FZ$, then we have $XPF=Q$, which is the equation from Lemma \ref{l:XPY=Q}. The first two cases correspond to two its solutions: either $X=\overline{F}ST$, or $F=ST\overline{X}$.

If $Y$ is shorter than $Z$, $Z=F'Y$, then we have $XP=QF'$, which is the equation from Lemma \ref{l:PX=YQ}, and the four last cases correspond to four its solutions. \hfill $\Box$

 \bigskip

\begin{lemma}\label{l:PXQ=R}
 Suppose that $PXQ=R$ for some palindromes $P$, $Q$, $R$. Then $X$ is a concatenation of (at most) three palindromes.
\end{lemma}
\noindent {\sc Proof.} 
This is a direct corollary of Theorem \ref{saar} applied twice: the palindromic length of $XQ$ is at most 2 and thus the palindromic length of $X$ is at most 3.
\hfill $\Box$

\begin{lemma}\label{l:P=QXRY}
Suppose that $P=QXRY$, where $P$, $Q$, $R$ are palindromes and $X$, $Y$ are arbitrary words. Then one of the following six cases holds: there exist three palindromes $S$, $T$ and $U$ and eventually an arbitrary word $Z$ such that
\begin{enumerate}
 \item either $X=S\overline{Y}TU$;
 \item or $X=ZST$ and $Y=\overline{Z}U$;
 \item or $X=SZT$ and $Y=U\overline{Z}$;
 \item or $X=ZS$ and $Y=T\overline{Z}U$;
 \item or $X=SZ$ and $Y=TU\overline{Z}$;
 \item or $Y=ST\overline{X}U$.
\end{enumerate}\end{lemma}
\noindent {\sc Proof.} Suppose first that $Y$ is shorter than a half of $Q$, that is, $Q=\overline{Y}Q'Y$ for a palindrome $Q'$. Then $Q'YXR$ is a palindrome, and we can use Lemma \ref{l:PXQ=R} to see that $YX$ is a concatenation of three palindromes. The six cases of the lemma statement correspond to six possible positions of the end of $Y$ in them: before the middle of the first one, after the middle in the first one, before the middle of the second one etc.

Suppose now that $Y$ is shorter than $Q$ but longer than its half: $Q=ZQ'\overline{Z}$ and $Y=Q'\overline{Z}$. Then $\overline{Z}XR$ is a palindrome, and we use Lemma \ref{l:P=QX} to see that $\overline{Z}X$ is a concatenation of two palindromes. Now we use Lemma \ref{l:XY=PQ} to get the following possibilities:
\begin{itemize}
 \item either $X=SZT$, and since $Y=Q'Z$, we get Case 3 of the statement;
 \item or $\overline{Z}=\overline{W}S$ and $X=WT$, and then $Y=Q'\overline{Z}=Q'\overline{W}S$, and after renaming variables, we get Case 4 of the statement;
 \item or $\overline{Z}=SW$ and $X=T'\overline{W}$, and then $Y=Q'\overline{Z}=Q'SW$, and after renaming variables, we observe Case 5;
 \item or $\overline{Z}=S\overline{X}T'$, and thus $Y=Q'S\overline{X}T'$, and this is Case 6.
\end{itemize}

At last, consider the possibility of $Y$ longer than $Q$. Then $Y=ZQ$, and $XRZ$ is a palindrome. So, we use Lemma \ref{l:XPY=Q} to see that either $X=\overline{Z}ST$ (and $Y=ZQ$), which falls into Case 2, or $Z=ST\overline{X}$ and $Y=ST\overline{X}Q$, which falls into Case 6. 

We have considered all the possibilities and got nothing except for the six listed cases. \hfill $\Box$

\begin{lemma}\label{l:PX=QYR}
 Suppose that $PX=QYR$ for some palindromes $P$, $Q$, $R$. Then there exist palindromes $S,T,U$ and eventually a word $Z$ that one of the following holds on the words $X$ and $Y$:
 \begin{enumerate}
  \item either $X=STYU$;
   \item or $X=SZT$ and $Y=UZ$;
  \item or $X=ZS$ and $Y=TUZ$;
  \item or $X=SZ$ and $Y=TZU$;
  \item or $Y=STXU$;
  \item or $X=STZ$ and $Y=ZU$.
 \end{enumerate}
\end{lemma}

\noindent {\sc Proof.} If $X$ is longer than $R$, $X=ZR$, then the equation is reduced to $PZ=QY$. Due to Lemma \ref{l:px=qy}, it means that either $Z=STY$ and thus $X=STYU$ (case 1), or $Y=STZ$ and $X=ZR$ (case 3 after renaming palindromes).

If $X$ is shorter than $R$ but longer than its half, we have $R=ZR'\overline{Z}$ and $X=R'\overline{Z}$ for a palindrome $R'$ and a word $Z$. So, the initial equation is reduced to $P=QYZ$, which means due to Lemma \ref{l:P=QX} that $YZ=ST$ for some palindromes $S$ and $T$. Now using Lemma \ref{l:XY=PQ}, we consider four possibilities, which, after renaming variables, give us cases 1, 6, 2 and 4 of the lemma.

If $X$ is shorter than a half of $R$, it means that $R=\overline{X}R'X$ for a palindrome $R'$. So, the equation is reduced to $P=QY\overline{X}R'$. Due to Lemma \ref{l:PXQ=R}, it means that $Y\overline{X}$ is a concatenation of three palindromes, and we have six cases to consider depending on the position of the boundary between $Y$ and $\overline{X}$ with respect to them. Five of them give cases which we have seen above, and the sixth corresponds to $Y=STXU$, which is case 5. \hfill $\Box$


\begin{lemma}\label{l:PX=QYRZ}
Suppose that $PX=QYRZ$, where $P$, $Q$, $R$ are palindromes and $X$, $Y$, $Z$ are arbitrary words. Then one of the following twelve cases holds: there exist three palindromes $S$, $T$ and $U$ and eventually an arbitrary word $C$ such that
\begin{enumerate}
 \item either $X=STYUZ$;
 \item or $X=SCTZ$, $Y=UC$;
 \item or $X=CSZ$, $Y=TUC$;
 \item or $X=SCZ$, $Y=TCU$;
 \item or $X=CZ$, $Y=STCU$;
 \item or $X=STCZ$, $Y=CU$;
 
 \item or $Y=SCTU$, $Z=\overline{C}X$;
 \item or $Y=CST$, $Z=\overline{C}UX$;
 \item or $Y=SCT$, $Z=U\overline{C}X$;
 \item or $Y=CS$, $Z=T\overline{C}UX$;
 \item or $Y=SC$, $Z=TU\overline{C}X$;
 \item or $Z=ST\overline{Y}UX$.
\end{enumerate}
\end{lemma}
\noindent {\sc Proof.} If $X$ is longer than $Z$, $X=WZ$, then the equation is immediately reduced to the equation $PW=QYR$ from Lemma \ref{l:PX=QYR}. The first six solutions correspond to six cases from that lemma, after renaming variables and adding the fact that $X=WZ$.

If $Z$ is longer than $X$, $Z=WX$, then we obtain the equation $P=QYRW$ from Lemma \ref{l:P=QXRY}. The last six cases correspond to its solutions after renaming variables and adding the fact that $Z=WX$. \hfill $\Box$

\begin{lemma}\label{l:PX=YQR}
 Suppose that $PX=YQR$ for some palindromes $P$, $Q$, $R$. Then there exist palindromes $S,T,U$ and eventually a word $Z$ that one of the following holds  on the words $X$ and $Y$:
 \begin{enumerate}
 \item either $X=SYTU$;
 \item or $X=ZST$, $Y=UZ$; 
 \item or $X=SZT$, $Y=ZU$; 
  \item or $X=ZS$, $Y=TZU$; 
   \item or $X=SZ$ and $Y=ZTU$; 
   \item or $Y=SXTU$. 
 \end{enumerate}
\end{lemma}
\noindent {\sc Proof.} If $Y$ is longer than $P$, $Y=PZ$, we immediately have $X=ZQP$, and this is Case 2. 

If $Y$ is shorter than $P$ but longer than its half, we have $P=CP'\overline{C}$ and $Y=CP'$ for some $C$. So, $\overline{C}X=QR$, and we have to consider the four cases from Lemma \ref{l:XY=PQ}. They give successively:

- either $Q=\overline{C}Q'C$, $X=Q'CR$, and this is Case 3; 

- or $Q=\overline{Z}Q'Z$, where $\overline{Z}Q'=\overline{C}$, and then $X=ZR$ and $Y=CP'=Q'ZP'$, which is Case 4; 

- or $\overline{C}=Q\overline{Z}$, where $R=\overline{Z}R'Z$ and then $X=R'Z$ and $Y=CP'=ZQP$, which is Case 5; 

- or $\overline{C}=Q\overline{X}R'$, where $R=\overline{X}R'X$, and thus $Y=CP'=RXQP'$, which is Case 6. 

At last, if $Y$ is shorter than a half of $P$, $P=YP'\overline{Y}$, then $P'\overline{Y}X=QR$, and due to Theorem \ref{saar}, it means that $\overline{Y}X$ is a concatenation of three palindromes. Checking all the six possible positions of the end of $\overline{Y}$ withing these three palindromes, from the shortest to the longuest $Y$, we get exactly the six cases in their order. \hfill $\Box$

\medskip

The next lemmas will be used in this paper for the case of the antimorphism $\varphi$ equal to the inverse. To give them more power, we state them in a general way. 

\begin{lemma}\label{l:APQAX=R}
Suppose that $APQ\varphi(A)X=R$ for some palindromes $P,Q,R$, a word $A$ and a bijective antimorphism $\varphi$ such that $\varphi(\varphi(a))=a$ for all $a\in \Sigma_g$. Then there exist palindromes $S,T,U$ and a word $Z$ such that one of the following holds:
\begin{enumerate}
 \item either $X=ZST\varphi(Z)U$;
 \item or $X=SZTU\varphi(Z)$.
\end{enumerate}
Note that these two representations are mirror images one of the other (after renaming variables).
\end{lemma}
\noindent {\sc Proof.} Suppose first that $A$ is shorter than $X$: $X=Y\overline{A}$ for a word $Y$. Then $PQ\varphi(A)Y$ is a central subpalindrome of $R$, and thus due to Theorem \ref{saar}, $\varphi(A)Y$ is a concatenation of three palindromes, say, $STU$. Here we have to consider the three cases for the position of the end of $\varphi(A)$ in them: is it in the first half of $S$, second half of $S$, first half of $T$, second half of $T$, first half of $U$, or second half of $U$.

All the cases are similar; we give a proof for the first of them. If $\varphi(A)$ is shorter than a half of $S$, we have $S=\varphi(A)S'\varphi(\overline{A})$ for a palindrome $S'$, and so $Y=S'\varphi(\overline{A})TU$. This gives $X=Y\overline{A}=S'\varphi(\overline{A})TU\overline{A}$, which is Case 2 of the statement.

All the other five cases are similar.

Now suppose that $A$ is longer than $X$: $A=\overline{X}B$ for a word $B$. Then $\varphi(A)=\varphi(B)\varphi(\overline{X})$, and we have $R=\overline{X}BPQ\varphi(B)\varphi(\overline{X})X$. So, the word $BPQ\varphi(B)\varphi(\overline{X})$ is a central subpalindrome $R'$ of $R$.

Applying $\varphi$ and the mirror image to the equation $BPQ\varphi(B)\varphi(\overline{X})=R'$, we get exactly the initial equation on $X$ with renamed constants. So, the form of its solutions stays the same. \hfill $\Box$

\begin{lemma}\label{l:APQAX=RY}
Suppose that $APQ\varphi(A)X=RY$ for some palindromes $P,Q,R$, a word $A$ and a bijective antimorphism $\varphi$ such that $\varphi(\varphi(a))=a$ for all $a\in \Sigma_g$. Then there exist palindromes $S,T,U$ and a word $Z$ such that one of the following holds:
\begin{enumerate}
 \item either $X=ZST\varphi(Z)UY$;
 \item or $X=SZTU\varphi(Z)Y$;
 \item or $Y=ZST\varphi(Z)UX$;
 \item or $Y=SZTU\varphi(Z)X$.
\end{enumerate}
\end{lemma}
\noindent {\sc Proof.} Suppose first that $X$ is longer than $Y$: $X=FY$. Then we have the equation $APQ\varphi(A)F=R$ which we have seen in Lemma \ref{l:APQAX=R} with renamed variables. So, either $F=ZST\varphi(Z)U$, and $X=FY$ corresponds to Case 1; or $F=SZTY\varphi(Z)$, and $X=FY$ corresponds to Case 2.

Now suppose that $Y$ is longer than $X$: $Y=FX$. Then the equation is reduced to 
\[APQ\varphi(A)=RF,\]
and we have two subcases.

- either $A$ is longer than $F$, which means that $\varphi(A)=CF$ for a word $C$ and thus $A=\varphi(F)\varphi(C)$. Then the equation can be reduced to
\[\varphi(F)\varphi(C)PQC=R,\]
and, applying to it the antimorphism $\varphi$, we get $\varphi(C)\varphi(Q)\varphi(P)CF=\varphi(R)$, which is the equation from Lemma \ref{l:APQAX=R} with renamed constants. So, either $F=ZST\varphi(Z)U$ and $Y=FX$ corresponds to Case 3; or $F=SZTU\varphi(Z)$, and for $Y=FX$ it means exactly Case 4.

- or $A$ is shorter than $F$, that is, $F=C\varphi(A)$ for a word $C$, and the equation is reduced to $APQ=RC$, which is exactly the equation from Lemma \ref{l:PX=YQR} with renamed variables and constants. So, we have to consider its six solutions.

Just for an example, let us consider here the solution 4: $C=BS$, $A=TBU$ for a word $B$ and palindromes $S, T, U$. Then $F=C\varphi(A)=BS\varphi(U)\varphi(B)\varphi(T)$, and since $\varphi$ of a palindrome is a palindrome, $Y=FX$ corresponds to the Case 3.

All the other solutions can be considered analogously and all give Cases 3 or 4, completing the proof. \hfill $\Box$

\section{Proof of Theorem \ref{t:2}}\label{s:t2}

First of all, clearly, if $P$ and $Q$ are palindromes and $A$ is a word, then $$APQA^{-1}=(AP\overline{A})(\overline{A}^{-1}QA^{-1}),$$ so, each word of the form indicated in the statement of the theorem is indeed a product of two palindromes. It remains to prove the opposite direction.

Let $U$ and $V$ be two palindromes. The reduced form of their product $UV$ is clearly obtained by the concatenation of a prefix of $U$ and a suffix of $V$. Both of them can be empty or equal to the initial palindrome, but what is important for us is if they are longer or shorter than a half of the initial palindrome. So, we have the following cases:

\medskip 
\noindent {\bf Case 1.} The remaining suffix of $V$ is longer than its half: it is equal to $V'Y$, where $V=\overline{Y}V'Y$ in the reduced form and $V'$ is a palindrome. Now there are two subcases:

\smallskip 
\noindent {\bf Case 1.1.} The remaining prefix of $U$ is longer that its half: it is equal to $XU'$, where $U=XU'\overline{X}$, and $U'$ is a palindrome. In this case, the annihilated words are $\overline{X}$ and $\overline{Y}$: $\overline{X}=(\overline{Y})^{-1}$, so, $Y=X^{-1}$, and the reduced form of $UV$ is $XU'V'Y=XU'V'X^{-1}$, confirming the theorem statement.

\smallskip 
\noindent {\bf Case 1.2.} The remaining prefix of $U$ is longer that its half: it is equal to $X$, where $U=XU'\overline{X}$, and $U'$ is a palindrome. So, the annihilated words are $U'\overline{X}$ and $\overline{Y}$; we have $U'\overline{X}=\overline{Y}^{-1}$, and thus $Y=(U')^{-1}X^{-1}$. The reduced form of $UV$ is $XV'Y=XV'(U')^{-1}X^{-1}$, confirming the theorem statement.

\medskip 
\noindent {\bf Case 2.}  The remaining suffix of $V$ is shorter than its half: it is equal to $Y$, where $V=\overline{Y}V'Y$ in the reduced form and $V'$ is a palindrome. Now there are also two subcases:

\smallskip 
\noindent {\bf Case 2.1.} The remaining prefix of $U$ is longer that its half: it is equal to $XU'$, where $U=XU'\overline{X}$, and $U'$ is a palindrome. This case is symmetric to the previous one.

\smallskip 
\noindent {\bf Case 2.2.} The remaining prefix of $U$ is longer that its half: it is equal to $X$, where $U=XU'\overline{X}$, and $U'$ is a palindrome. So, the annihilated words are $U'\overline{X}$ and $\overline{Y}V'$, and we have $U'\overline{X}=(\overline{Y}V')^{-1}=(V')^{-1}\overline{Y}^{-1}$ (in reduced forms). Applying Lemma \ref{l:px=qy}, we get $\overline{X}=ST\overline{Y}^{-1}$ or $\overline{Y}^{-1}=ST\overline{X}$, that is, $X=Y^{-1}TS$ or $Y^{-1}=XTS$ and thus $Y=S^{-1}T^{-1}X^{-1}$ for some palindromes $S$ and $T$.

The reduced form of $UV$ is $XY$; according to the result of the previous paragraph, $XY$ is equal to $Y^{-1}TSY$ or to $XS^{-1}T^{-1}X^{-1}$ (in the reduced form) for some palindromes $S$ and $T$. 

We have considered all the possible cases and proved Theorem \ref{t:2}. \hfill $\Box$

\section{Proof of Theorem \ref{t:3}}\label{s:t3}

This proof follows the method of the previous one, but is much more technical. Everywhere in the proof, the letters $P, Q, R, S, T, U, V, W$ with indices are reserved for palindromes, and all other capital latin letters mean general words, by default in reduced form.

From now on, the two possible reduced forms, $ABPQB^{-1}R \overline{A}$ and $APBQRB^{-1}\overline{A}$ will be referred to as the $ABP-$word and thr $APB-$word. We will obtain each of them many times with renamed variables and will recognize them like that.

First of all, 
$$ABPQB^{-1}R \overline{A}=(ABP\overline{B}\;\overline{A})(\overline{A}^{-1}\overline{B}^{-1}QB^{-1}A^{-1})(AR\overline{A})$$ 
and 
$$ APBQRB^{-1}\overline{A}=(AP\overline{A})(\overline{A}^{-1}BQ\overline{B}A^{-1})(A\overline{B}^{-1}RB^{-1}\overline{A}),$$
so that in both cases, the initial word is a product of three palindromes. It remains to prove the opposite direction: every product of three palindromes has one of the two reduced forms above.

Let $U$, $V$, $W$ be three palindromes in $\langle\Sigma\rangle$. We study their product $UVW$. Due to Theorem \ref{t:2}, the reduced form of $UV$ is $APQA^{-1}$ for some palindromes $P$ and $Q$ and a word $A$. It remains to study all the possible cases of what prefix $Y^{-1}$ of $W$ annihilates with what suffix $Y$ of $APQA^{-1}=XY$.

The main two cases to consider are: Is that annihilated prefix longer or shorter than a half of $W$?

\bigskip
\noindent {\bf Case 1}: The prefix is short, namely, the reduced form of $W$ is $W=Y^{-1}W'\overline{Y}^{-1}$ for a palindrome $W'$ (here $C$ and $W'$ may be empty).

In this case, the reduced form of $UVW$ is  $XW'\overline{Y}^{-1}$ for some decomposition $APQA^{-1}=XY$. It remains to list all possible situations for $X$ and $Y$.

\medskip
\noindent {\bf Case 1.1}: $Y$ is a suffix of $A^{-1}$, that is, $A^{-1}=A'Y$ and thus $X=APQA'$. Here $A=Y^{-1}(A')^{-1}$, and thus $XW'\overline{Y}^{-1}=Y^{-1}(A')^{-1}PQA'W'\overline{Y}^{-1}$, which is an $ABP-$word.

\medskip 
\noindent {\bf Case 1.2}: $Y$ is a suffix of $QA^{-1}$ which does not contain the middle of $Q$ but touches $Q$, so that $Q=CQ_1\overline{C}$ and $Y=\overline{C}A^{-1}$. Then $X=APCQ_1$, $\overline{Y}^{-1}=C^{-1}\overline{A}$, and thus the reduced form of $UVW$ is $XW'\overline{Y}^{-1}=APCQ_1W'C^{-1}\overline{A}$, which is an $APB-$word.

\medskip 
\noindent {\bf Case 1.3}: $Y$ is a suffix of $QA^{-1}$ which does contain the middle of $Q$, so that $Q=CQ_1\overline{C}$ and $Y=Q_1\overline{C}A^{-1}$. Then $X=APC$, $\overline{Y}^{-1}=Q_1^{-1}C^{-1}\overline{A}$, and thus the reduced form of $UVW$ is $XW'\overline{Y}^{-1}=APCW'Q_1^{-1}C^{-1}\overline{A}$, which is an $APB-$word.

\medskip 
\noindent {\bf Case 1.4}: $Y$ is a suffix of $PQA^{-1}$ which does not contain the middle of $P$ but touches $P$, so that $P=CP_1\overline{C}$ and $Y=\overline{C}QA^{-1}$. Then $X=ACP_1$, $\overline{Y}^{-1}=C^{-1}Q^{-1}\overline{A}$, and thus the reduced form of $UVW$ is $XW'\overline{Y}^{-1}=ACP_1W' C^{-1}Q^{-1}\overline{A}$, which is an $ABP-$word.

\medskip 
\noindent {\bf Case 1.5}: $Y$ is a suffix of $PQA^{-1}$ which does contain the middle of $P$, so that $P=CP_1\overline{C}$ and $Y=P_1\overline{C}QA^{-1}$. Then $X=AC$, $\overline{Y}^{-1}=P_1^{-1}C^{-1}Q^{-1}\overline{A}$, and thus the reduced form of $UVW$ is $XW'\overline{Y}^{-1}=ACW'P_1^{-1}C^{-1}Q^{-1}\overline{A}$, which is an $ABP-$word.

\medskip 
\noindent {\bf Case 1.6}: $Y$ is a suffix of $APQA^{-1}$ touching $A$, that is, $A=XC$ and 
$Y=CPQA^{-1}$. Then $\overline{Y}^{-1}=\overline{C}^{-1}P^{-1}Q^{-1}\overline{A}=\overline{C}^{-1}P^{-1}Q^{-1}\overline{C}\;\overline{X}$, and the reduced form of $UVW$ is $XR\overline{C}^{-1}P^{-1}Q^{-1}\overline{C}\;\overline{X}$ , which is a $APB-$word.
\medskip
The case of a short annihilated prefix of $W$ is over; now let us treat the more complicated one.

\bigskip
\noindent {\bf Case 2}: After multiplication to $UV=APQA^{-1}=XY$, more than a half of the length of $W$ annihilates with it: $W=CW'\overline{C}$ for a palindrome $W'$, and  $CW'$ is equal to $Y^{-1}$, whereas the reduced form of $UVW$ is $X\overline{C}$. Now for all possible decompositions $APQA^{-1}=XY$ we have to solve the equation 
$$Y^{-1}=CW'$$ 
and to express $X\overline{C}$.

\medskip
\noindent {\bf Case 2.1}: $Y$ is a suffix of $A^{-1}$: $A^{-1}=BY$ and thus $A=Y^{-1}B^{-1}$. Since $Y^{-1}=CW'$, it means that $X=APQB=CW'B^{-1}PQB$ and thus $X\overline{C}= CW'B^{-1}PQB\overline{C}$, which is a $APB-$word.

\medskip
\noindent {\bf Case 2.2}: $Y$ is a suffix of $QA^{-1}$ which touches $Q$ but not its middle: $Y=BA^{-1}$, where $Q=\overline{B}Q'B$ for some middle palindrome $Q'$. So, $Y^{-1}=AB^{-1}=CW'$, and our goal is to express $X\overline{C}=AP\overline{B}Q'\overline{C}$.

The equation $AB^{-1}=CW'$ is that from Lemma \ref{l:xy=zp} with renamed variables. So, according to it, we have the three possible subcases:

\begin{itemize}
 \item either $C=AD$ and $B^{-1}=DW'$ for some $D$; then $X\overline{C}=AP\overline{B}Q'\overline{C}=AP\overline{D}^{-1}(W')^{-1}Q'\overline{D}\;\overline{A}$, which is an $APB-$word;
 \item or $A=CD$ and $B^{-1}=W''\overline{D}$ for a word $D$ and a palindrome $W''$; then $X\overline{C}=AP\overline{B}Q'\overline{C}= CDP(W'')^{-1}D^{-1}Q'\overline{C}$, which is a $ABP-$word;
 \item or $A=C\overline{B}^{-1}W''$ for a palindrome $W''$, and then
 $X\overline{C}=AP\overline{B}Q'\overline{C}=C\overline{B}^{-1}W''P\overline{B}Q'\overline{C}$, which is an $ABP-$word.
\end{itemize}

\medskip
\noindent {\bf Case 2.3}: $Y$ is a suffix of $QA^{-1}$ which does touch the middle of the palindrome $Q$: $Y=Q_1DA^{-1}$, where $Q=\overline{D}Q_1D$, and at the same time, $Y=(CW')^{-1}=(W')^{-1}C^{-1}$, where all the forms are reduced. Our goal here is to express $X\overline{C}=AP\overline{D}\overline{C}$.

The equation $Q_1DA^{-1}=(W')^{-1}C^{-1}$, where $Q_1$ and $W'$ are palindromes and the other variables are not obliged to be so, is in fact the equation from Lemma \ref{PXY=QZ} with renamed variables. 

So, we have the six cases to consider, exactly as above.

For example, the first case gives $C^{-1}=STDA^{-1}$, which means $\overline{C}=S^{-1}T^{-1}\overline{D}^{-1}\overline{A}$ and thus

 \[AP\overline{D}\;\overline{C}=AP\overline{D}S^{-1}T^{-1}\overline{D}^{-1}\overline{A},\]
 which is a $APB-$word.
 
 The other cases are similar.

 \medskip
\noindent {\bf Case 2.4}: $Y$ is a suffix of $PQA^{-1}$ which touches $P$ but not its center: $P=\overline{D}P_1D$ and $Y=DQA^{-1}$. As always, at the same time, $Y=(W')^{-1}C^{-1}$, and we have to consider the following equation on reduced forms:
\[DQA^{-1}=(W')^{-1}C^{-1}.\]
Our goal here is to express $X\overline{C}=A\overline{D}P_1\overline{C}$.

The equation above is the equation from Lemma \ref{l:XPY=QZ}, and it is sufficient to consider its 6 solutions listed in the lemma. 

For example, the first one after renaming variables is $D=\overline{F}ST$ and $A^{-1}=FC^{-1}$, that is, $\overline{D}=TSF$ and $A=C F^{-1}$. So, $A\overline{D}P_1\overline{C}=C F^{-1}TSFP_1\overline{C}$ for some palindromes $S$, $T$ and a word $F$, which is an $ABP-$word.


The other five cases are considered analogously.

 \medskip

\noindent {\bf Case 2.5}: $Y$ is a suffix of $PQA^{-1}$ which touches the center of $P$: $P=DP_1\overline{D}$ and $Y=P_1\overline{D}QA^{-1}$. So, the equation to solve here is
\[P_1\overline{D}QA^{-1}=(W')^{-1}C^{-1},\]
which is the equation from Lemma \ref{l:PX=QYRZ} with renamed variables.
Our goal is to express $A D \overline{C}$, so, we have to consider the twelve cases separately. 

Just for an example,  Case 10 of the lemma means that $\overline{D}=FS$ and $A^{-1}=T\overline{F}UC^{-1}$ for some palindromes $S$, $T$, $U$ and a word $F$. So, $D=S\overline{F}$, $A=CU^{-1}\overline{F}^{-1}T^{-1}$, and 
\[AD\overline{C}=CU^{-1}\overline{F}^{-1}T^{-1}S\overline{F}\overline{C},\]
which is a $APB-$word. The other 11 cases are considered analogously.

\medskip
\noindent {\bf Case 2.6}: $Y$ is a suffix of $APQA^{-1}$ which touches the first occurrence of $A$: $A=FG$, $A^{-1}=G^{-1}F^{-1}$, and $Y=GPQG^{-1}F^{-1}$. So, here we solve the equation 
\[GPQG^{-1}F^{-1}=(W')^{-1}C^{-1},\]
and our goal is to express $F\overline{C}$.

The equation is that from Lemma \ref{l:APQAX=RY} with renamed variables and the inverse as $\varphi$, so, we have four cases to consider. For example, the first of them gives

\[F^{-1}=ZSTZ^{-1}UC^{-1},\]
that is, $F=CU^{-1}ZT^{-1}S^{-1}Z$, and thus $F\overline{C}=CU^{-1}ZT^{-1}S^{-1}Z\overline{C}$ is an $APB$-word.

The other four cases are analogous, and the proof of the theorem is complete. \hfill $\Box$

\section{Conjectures}
The characterisations of words of group palindromic complexity $2$ and $3$ obtained in this paper suggest that checking if this complexity exceeds $k$ or not can be decidable for every $k$, and that the decision procedure involves solving semigroup equations with palindromes and antimorphisms.

 Moreover, we can easily imaging how the next characterizations look like. Indeed, the results of Theorems \ref{t:2} and \ref{t:3} mean exactly that all the possible reduced forms look as if the centers of initial palindromes are not touched by reduction. As we have shown above,  
$$APQA^{-1}=(AP\overline{A})(\overline{A}^{-1}QA^{-1}),$$
$$ABPQB^{-1}R \overline{A}=(ABP\overline{B}\;\overline{A})(\overline{A}^{-1}\overline{B}^{-1}QB^{-1}A^{-1})(AR\overline{A})$$ 
and 
$$ APBQRB^{-1}\overline{A}=(AP\overline{A})(\overline{A}^{-1}BQ\overline{B}A^{-1})(A\overline{B}^{-1}RB^{-1}\overline{A}),$$
so that $P$, $Q$, $R$ are just centers of the initial palindromes. The difference between the two last forms is in which of the two reductions is longer, the first or the second one.

We conjecture that the same structure holds in words decomposable to any number $k$ of palindromes: for every allowed reduced form, there exist a decomposition to $k$ palindromes where the centers of the initial palindromes did not change with reduction. So, to find all the possible forms, it is sufficient to consider all the $(k-1)!$ possible orderings of lengths reduced between every two palindromes.

In particular, it would mean that a word is decomposable to $4$ palindromes only if it has one of the following $5$ forms:
\begin{eqnarray*}
&& APBQCRSC^{-1}\overline{B}A^{-1}; \\
&& APBCQRC^{-1}S\overline{B}A^{-1};\\
&& ABPQB^{-1}CRSC^{-1}A^{-1};\\
&& ABPCQRC^{-1}\overline{B}SA^{-1}\\
&& ABCPQC^{-1}R\overline{B}SA^{-1}.
\end{eqnarray*}
Note that the central form, structurally different from the others, correspond to both cases when the shortest reduction is in the middle, no matter which of the side ones is longer.

This fact and the first values $1,1,2$ and perhaps $5$ of the number of possible expressions for $k=1,2,3,4$, suggest that perhaps the reduced forms may somehow be related to Dick paths, and their numbers are Catalan numbers.

Regardless of whether or not it is true, it is clear that both the number of possible forms and the number of cases to consider for sthe respective theorem grow extremely fast with $k$. So, our method cannot be directly used to find the reductions for every $k$. 

However, in some sense, in this paper we reduce decidability of the problem ``Is a word decomposable to $k$ group palindromes'' to solving semigroup equations with palindromes and antimorphisms. If someone manages to extend to such equations the known methods, like the Makanin algorithm, it will mean that the problem is decidable.

\end{document}